\documentclass{article}

\usepackage{amsmath}
\usepackage{amssymb}
\usepackage{latexsym}
\usepackage{theorem}
\usepackage{epsfig}
\usepackage{subfigure}
\usepackage{varioref}

\newtheorem{theorem}{Theorem}[section] 
\newtheorem{lemma}[theorem]{Lemma}

\newtheorem{corollary}[theorem]{Corollary}
\newtheorem{conjecture}[theorem]{Conjecture}

{\theorembodyfont{\rmfamily}
\theoremstyle{plain}
\newtheorem{definition}[theorem]{Definition}
\newtheorem{example}[theorem]{Example}

\newtheorem{ack}{Acknowledgements}
}

\newcommand{\A}{\ensuremath{\mathcal A}}
\newcommand{\B}{\ensuremath{\mathcal B}}
\newcommand{\codim}{\operatorname{codim}}
\newcommand{\I}{\ensuremath{\mathcal I}}
\newcommand{\E}{\ensuremath{\mathcal E}}
\newcommand{\J}{\ensuremath{\mathcal J}}
\newcommand{\C}{\ensuremath{\mathbb{C}}}
\renewcommand{\P}{\ensuremath{\mathbb{P}}}
\newcommand{\Z}{\ensuremath{\mathbb{Z}}}

\newcommand{\K}{\ensuremath{\mathbb{K}}}

\newcommand{\cl}{\ensuremath{c\ell}}
\newcommand{\lc}{\ensuremath{\ell c}}
\renewcommand{\to}{\ensuremath{\longrightarrow}}
\newcommand{\al}{\ensuremath{\alpha}}
\newcommand{\we}{\ensuremath{\wedge}} 
\renewcommand{\bar}[1]{\ensuremath{\overline{#1}}}
\newcommand{\rk}{\operatorname{rk}}
\newcommand{\nb}{\operatorname{\it nbb}}
\newcommand{\NB}{\operatorname{\it NBB}}
\newcommand{\nc}{\operatorname{\it nbc}}
\newcommand{\nbb}{{\em nbb}}
\newcommand{\nbc}{{\em nbc}}
\newcommand{\NBB}{{\em NBB}}
\newcommand{\qed}{\hfill \mbox{$\Box$}\medskip}
\newenvironment{proof}{\noindent {\it proof:}}{\qed \par}

\renewcommand{\theenumi}{\rm \theenumi}

\renewcommand{\theenumi}{\roman{enumi}}

\title{Line-closed matroids, quadratic 
algebras, and formal arrangments} 

\author{Michael Falk} 
\begin{document}
\maketitle

\begin{abstract}

Let $G$ be a matroid on ground set \A. The Orlik-Solomon algebra 
$A(G)$ is the quotient of the exterior algebra \E\ on \A\ by the ideal 
\I\ 
generated by circuit boundaries. The quadratic closure $\bar{A}(G)$
of $A(G)$ is the quotient of \E\ by the 
ideal generated by the degree-two component of \I.
We introduce the notion of \nbb\ set in $G$, determined by a linear 
order on \A, and show that the 
corresponding monomials are linearly independent in the quadratic 
closure $\bar{A}(G)$. As a consequence, $A(G)$ is a quadratic algebra 
only if $G$ is line-closed. An example of S.~Yuzvinsky proves the 
converse false. These results generalize to the degree $r$ closure of 
$\A(G)$.

The motivation for studying line-closed matroids grew out of the 
study 
of formal arrangements. This is a geometric 
condition necessary for \A\ to be free and for the complement $M$ of 
\A\ to be a $K(\pi,1)$ space.  Formality of \A\ is also necessary for 
$A(G)$ to be a quadratic algebra. We clarify the relationship between 
formality, 
line-closure, and other matroidal conditions related to formality.  We give 
examples to show that line-closure of $G$ is not necessary or sufficient for 
$M$ 
to be a $K(\pi,1)$, or for \A\ to be free.
\end{abstract}

\begin{section}{Introduction}
\label{intro}
Let \K\ be a field. An  {\em arrangement} is a finite set $\A=\{H_1, 
\ldots, 
H_n\}$ of linear 
hyperplanes in $V=\K^\ell$. Each $H_i$ is the kernel of a linear form 
$\al_i : V \to \K$, unique up to nonzero scalar multiple. 
Let $[n]$ denote the set $\{1,\ldots,n\}$ and $2^{[n]}$ the 
set of subsets of 
$[n]$.

A coordinate-free combinatorial model of the arrangement \A\ is 
provided by the 
{\em underlying matroid} of \A, which we denote by $G(\A)$, or
simply $G$. This matroid 
contains 
the same information as the intersection lattice $L(\A)$ -- see
\cite{OT}. 
By definition the matroid $G$ is the collection of dependent subsets 
of the set 
of defining forms $\{\al_1, \ldots, \al_n\}$. We identify these 
subsets 
with the corresponding sets of labels. Then it is easy to see that 
$$G=\{ \, S 
\subseteq [n] \ | \ \codim(\bigcap_{i\in S} H_i)<|S| \, \}.$$
Elements of $G$ are called {\em dependent sets}, and elements of 
$2^{[n]} - G$ 
are {\em independent sets}.  The projective point configuration $\A^* 
\subseteq 
\P(V^*)$ determined by $\{\al_1, \ldots, \al_n\}$ is called a {\em 
projective 
realization} of $G$. 

There are several other data besides the dependent sets which suffice 
to 
determine $G$ uniquely. Among these are the {\em 
circuits} 
of $G$, which are the minimal dependent sets, and the {\em bases} of 
$G$, which 
are the maximal independent sets. 
Besides these, we single out two functions which also uniquely 
determine $G$. 
The {\em rank function} $\rk: 2^{[n]} \to \Z$  is given by 
$\rk(X)=\codim(\bigcap_{i\in X} H_i)$. In the abstract setting, 
$\rk(X)$ is the 
(unique) size of a maximal independent subset of $X$. The rank 
$\rk(G)$ of $G$ 
is $\rk([n])$.
The {\em closure operator} $\cl: 2^{[n]} \to 2^{[n]}$ , given by 
$$\cl(X)=\{ 
i\in [n] \ | \ \rk(X\cup\{i\})=\rk(X)\},$$ also uniquely determines 
$G$. 

We refer the reader to the recent survey \cite{FR4} for more 
discussion of the role of matroid theory in the study of complex hyperplane 
arrangements.

A set $S$ is {\em closed} if $\cl(S)=S$. 
Closed sets are also
called {\em flats}. A flat corresponds to the collection of 
hyperplanes in \A\ 
containing a fixed subspace of $\K^\ell$, or equivalently, the 
intersection of 
the 
point configuration $\A^*$ with a fixed projective subspace of 
$\P(V^*)$. 
The set of flats, ordered by inclusion, forms a geometric lattice  
$L(G)$
isomorphic to the intersection lattice $L(\A)$. The flats of rank one are the 
singletons, called {\em points}. Flats of rank two are called {\em lines}. 
This 
terminology is 
natural with regard to the dual projective point configuration 
$\A^*.$  

Let $\K=\C$. The {\em complement} of \A\ is $V - \bigcup_{i=1}^n H_i,$ denoted 
by 
$M$.  The cohomology $H^*(M)$ is isomorphic to 
the {\em Orlik-Solomon algebra} $A(G)$ of the underlying matroid $G$, defined 
in 
the next section.  
Study of the lower central series of $\pi_1(M)$ \cite{F4,Y5} leads to the 
consideration of arrangements for which the cohomology algebra $H^*(M)$, or 
equivalently, the 
Orlik-Solomon algebra $A(G)$, is quadratic. Here $A(G)$ is {\em quadratic} if 
it 
has a presentation in which 
all relations have degree two.   
While this condition depends only on $G$,
the underlying combinatorial meaning has never been 
understood.  

The best results in this direction involve the notion of formality. An 
arrangement is {\em formal} if it is uniquely determined up to linear 
isomorphism 
by the dependence {\em relations} yielding dependent sets of rank two in $G$. 
This is a geometric, non-matroidal condition, and is a necessary condition for 
$A(G)$ to be quadratic. 
In looking for a matroidal analogue of formality we were led naturally to the 
study of line-closed matroids. 

\begin{definition}
A subset $S\subseteq [n]$ is {\em line-closed} if $\cl(\{i,j\}) 
\subseteq S$ for 
every $i,j \in S$. The matroid $G$ is {\em line-closed} if every 
line-closed set 
is closed.
\end{definition}

In attempting to sort out how this property 
fits in with other properties related to formality, we were led to 
the following.
\begin{conjecture} 
$G$ is line-closed  if and only if $A(G)$ is quadratic. 
\label{boston}
\end{conjecture}
In this paper 
we prove half of this conjecture, that $A(G)$ quadratic implies $G$ 
line-closed, for arbitrary coefficient fields \K.

The author sketched this 
proof and stated Conjecture~\ref{boston} in a lecture at the workshop 
``Arrangements in Boston" in 1999 \cite{F13}. Subsequently S.~Yuzvinsky 
found a counter-example for the full conjecture (at least for $\K=\C$), a 
line-closed 
matroid 
whose Orlik-Solomon  algebra is not quadratic. We exhibit Yuzvinsky's 
example, and refer the reader to the companion paper \cite{DY} for 
details on Yuzvinsky's approach, and for a stronger 
condition, also necessary but not sufficient for quadraticity of 
$A(G)$.

\medskip
One can define a quadratic algebra $\bar{A}(G)$ and a surjection 
$\bar{A}(G) \to A(G)$ which is an isomorphism if and only if $A(G)$ 
is 
quadratic; $\bar{A}(G)$ is called the  {\em quadratic closure} of 
$A(G)$.
Our main theorem follows from a more general 
construction, a partial generalization of the well-known 
\nbc\ (``no broken circuits'') basis for $A(G)$. We generalize one of 
the characterizations of \nbc\ sets to the lattice of line-closed 
sets 
of $G$. The result is the notion of \nbb\ {\em set.} 
We show that the monomials corresponding to \nbb\ sets are 
linearly independent in the quadratic closure $\bar{A}(G)$. 
In contrast to the situation for \nbc\ sets, the number 
of \nbb\ sets is not independent of the linear ordering of the atoms.
But the collection of \nbb\ sets will include all of the \nbc\ sets, 
for any given linear ordering. The two collections coincide, for 
every 
linear ordering, if and only if $G$ is line-closed. Thus, if $G$ is 
not line-closed, $\bar{A}(G)$ must be strictly bigger than $A(G)$, 
so $A(G)$ is not quadratic. The entire development generalizes to any 
degree, with the line-closed sets and quadratic 
closure replaced by $r$-closed sets and degree $r$ closure. The main 
theorem and its generalization are developed and proved in Section~\ref{main}.

The problem of finding a (monomial, or 
combinatorial) basis for $\bar{A}(G)$ is an interesting problem with some 
applications to lower central series calculations. Yuzvinsky's example shows 
that 
there may be no linear ordering for which the \nbb\ monomials 
form a basis of $\bar{A}(G)$.
Our definition of \nbb\ set
is a special case of the \NBB\ (``no bounded below'') 
sets of A.~Blass and B.~Sagan \cite{BSag}, for the lattice of line-closed sets 
of 
$G$, with a linear ordering on the set of 
atoms. Blass and Sagan define \NBB\ sets for finite atomic lattices 
with an arbitrary partial order on the set of atoms. Although general \NBB\ 
monomials for the lattice of line-closed sets are 
not linearly independent in $\bar{A}(G)$, we present a partial generalization 
of 
our main result to non-linear orderings, possibly yielding better 
lower bounds on $\dim 
\bar{A}(G)$ for $G$ of rank four or greater. 

\medskip
Much of the research in complex hyperplane arrangements focuses on the 
extent to which properties of the complement $M$ as a topological space or 
algebraic variety are determined by
the combinatorial structure of $G$.  In particular, two important open 
problems 
are whether asphericity of $M$ \cite{FR1} or freeness of \A\ \cite{OT} are 
dependent only on $G$ -- see \cite{FR1} and \cite{OT}. Formality of \A\ is 
also 
a necessary condition for each of 
these two properties. Thus  attempts 
were made to replace the definition of formality with 
some stronger purely 
combinatorial notion -- line-closure is one example. In the last section we 
give 
several other natural candidates for combinatorial analogues of formality.
each of them is stronger than formality. 
We establish the relationships among these various notions and show by 
example that in fact none of them have true topological implications. The 
discussion leads to an interesting conjecture concerning matroids which are 
determined by their points and lines.
\end{section}

\begin{section}{Quadratic closure and \nbb\ sets}
\label{main}
We will use the matroid-theoretic terminology developed in the 
introduction 
without further comment. The reader is referred to \cite{Wh1,Ox} for 
further 
background.

We begin with the definition of the Orlik-Solomon algebra 
$A(G)$ of a matroid $G$ on ground set $[n]$. For the remainder of the paper, 
let 
\K\ be any field, or indeed any commutative ring. Assume $G$ has no loops or 
multiple points.

Let $\E=\Lambda(V)$, the 
exterior algebra generated by $1$ and  $\{e_i \ | \ 1\leq i \leq 
n\}$, with the usual grading by degree. If $S=(i_1,\ldots,i_p)$ is an 
ordered $p$-tuple we denote the  product
$e_{i_1}\cdots e_{i_p}$ by $e_S$. We occasionally use the same 
notation when $S$ is an unordered set -- in this case $e_S$ is 
well-defined up to sign.

Define the linear mapping $\partial : \E^p \to \E^{p-1}$ by 
$$\partial(e_{i_1} \cdots
e_{i_p})=
\sum_{k=1}^p
(-1)^{k-1}e_{i_1} \cdots \widehat{e}_{i_k}\cdots e_{i_p},$$
where $\widehat{\hspace{.5em}}$ indicates an omitted factor. Then $\partial$M is 
a graded derivation, that is, $$\partial (x\we y)=\partial x\we 
y+(-1)^{\deg(x)}x\we\partial y$$ for homogeneous $x,y\in \E$.

Let $\I$ denote the ideal of \E\ generated by $\{\partial e_S \ | \ S 
\ \text{is dependent}\}$. 
\begin{definition} The {\em Orlik-Solomon algebra} $A=A(G)$ of $G$ is 
the quotient $\E/\I$.
\end{definition}
Since \I\ is generated by homogeneous elements, both \I\ and $A$ 
inherit gradings from \E. We denote the image of $e_S$ in $A$ by 
$a_S$.
The topological significance of $A$ is 
given in the following.

\begin{theorem}[{\rm \cite{OS1}}] If \A\ is an arrangement in 
$\C^\ell$ 
with 
complement $M$ and underlying matroid $G$, then $A(G)\cong H^*(M,\K)$.
\label{orsol}
\end{theorem}

\begin{paragraph}{The quadratic closure of $A(G)$}

\begin{definition} A graded algebra $U$ is {\em 
quadratic} if $U$ has a presentation with generators of degree one and 
relations of degree at most two.
\end{definition}

Let \J\ denote the ideal of \E\ generated by $\I^2$, the degree two 
part of the relation ideal \I. Because \E\ itself is quadratic, the 
Orlik-Solomon 
algebra $A$ will be quadratic if and only if $\J=\I$.
More generally the quotient $\E/\J$ is called the {\em quadratic 
closure} 
of $A(G)$, denoted $\bar{A}(G)$, or sometimes $\bar{A}$.

Quadratic Orlik-Solomon algebras appear in the study of complex 
arrangements, in the rational homotopy theory of the complement $M$ \cite{F4}
and the Koszul property of $A(G)$ \cite{Y5}. 
Rational homotopy theory 
provides a connection with the lower central series of the 
fundamental 
group. In that vein an invariant $\phi_3$ of $A(G)$ was introduced in 
\cite{F6}, defined as follows:
$$\phi_3(G)=\mbox{nullity}\left(\delta: \E^1\times \I^2 \to 
\E^3\right),$$
where $\delta$ is multiplication in $\E$. When $\K=\C$ and $G$ is the 
matroid of an arrangement with complement $M$, $\phi_3(G)$ is the 
rank of the third factor in lower central series of $\pi_1(M)$. 
Because the image of $\delta$ is precisely $\J^3$, the cokernel of 
$\delta$ is $\bar{A}\,^3$, and we find a simple relationship between 
$\phi_3(G)$ and $\dim(\bar{A}\,^3)$. The proof of the identity is left 
as 
an exercise.

\begin{theorem} 
\begin{equation*}
\begin{split}
\phi_3(G)&=\dim(\bar{A}\,^3) + n\dim(\I^2)-{n\choose 3}\\
&=2{{n+1}\choose 3} - n\dim(A^2) 
+\dim(\bar{A}\,^3)
\end{split}
\end{equation*}
\end{theorem}

The preceding identity can be stated in a simpler way, 
indicating that 
$\dim(\bar{A}\,^3/A^3)$ measures of the failure of the 
LCS (lower central series) formula relating the ranks of lower 
central series factors of $\pi_1(M)$ to the betti numbers 
$\dim(A^p)$, $p\geq 0$. 
Let $\gamma_3=\dim(A^3) +n\dim(\I^2)-{n\choose 3}$. Then, according to 
\cite{F8}, 
$\gamma_3$ is the value of $\phi_3$ predicted by the LCS formula for the given 
values of $\dim(A^p)$, and $\phi_3\geq \gamma_3$ with equality if and only if 
$\I^3=\J^3$. 

\begin{corollary} $\dim(\bar{A}\,^3) - \dim(A^3)=\phi_3-\gamma_3$
\label{count}
\end{corollary}

\end{paragraph}

\begin{paragraph}{The line-closure of a matroid}
Next we refine further the material on line closure from Section~\ref{intro}. 
Let us define an idempotent, order-preserving closure
operator on subsets of $S\subseteq [n]$ using line-closed sets: the 
{\em line-closure} $\lc(S)$ is by definition the intersection of the 
line-closed sets containing $S$. Since closed sets are automatically 
line-closed, $\lc(S)\subseteq \cl(S)$. We will consider the 
combinatorial 
structure consisting of the set $[n]$ equipped with the closure 
operator $\lc: 2^{[n]} \to 2^{[n]}$ to be the {\em line-closure} of 
the matroid $G$, and denote it by $\bar{G}$. This set system $\bar{G}$ 
will not be a 
matroid in general, because the operator 
$\lc$ fails to satisfy the Steinitz exchange axiom. The arguments and 
constructions in this section are seriously affected by this defect.  
Clearly $G$ is 
line-closed if and only if $\bar{G}=G$.

The collection of line-closed sets, partially-ordered by inclusion, 
will be denoted by $\bar{L}(G)$.  We call $\bar{L}(G)$ the {\em 
line-closure} of $L(G)$.  The poset $\bar{L}(G)$ is a lattice 
\cite[Section 2]{Rota}
in which every 
element is a join of atoms. But $\bar{L}(G)$ is not a graded 
lattice, as 
the example below shows.  This is a reflection of the failure of the 
exchange axiom. Again, $G$ is line-closed if and only if 
$\bar{L}(G)=L(G)$.

\begin{example} Let $G$ be the rank-three matroid on $[6]$ with 
rank-two circuits $\{1,2,3\}, \{3,4,5\}$, and $\{1,5,6\}$.  This 
example, pictured in Figure~\vref{3wheel}, is the ``rank-three wheel" 
\cite{Ox}.

\begin{figure}

\begin{center}
\epsfig{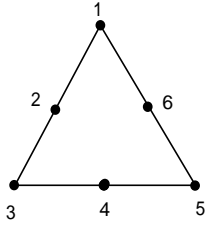}
\caption{The rank-three wheel}
\label{3wheel}
\end{center}
\end{figure}

Then there are two maximal chains in $\bar{L}(G)$ of different 
lengths, 
namely $$\emptyset<1<123<123456,
\ \text{and}$$
$$\emptyset<2<24<246<123456.$$
\label{wheel}
\qed
\end{example}

\end{paragraph}

\begin{paragraph}{\nbc\ and \nbb\ sets}
Fix a linear order of the ground set $[n]$.  A broken circuit of $G$ 
is a set of the form $C-\min(C)$, where
$C$ is a circuit of $G$.  An \nbc\ set of $G$ is a subset of $[n]$ 
which contains no broken circuits.  The collection of \nbc\ sets of 
$G$ 
will be denoted $\nc(G)$. The set of elements of $\nc(G)$ of 
cardinality $p$ is denoted $\nc^p(G)$. The dependence on the linear 
order of $[n]$ is suppressed in the notation. For a flat $X$ of $G$, 
let 
$\nc_X(G)$ denote the set of \nbc\ sets with closure equal to $X$. 

Among the properties of $\nc(G)$ we highlight the following. For proofs and a 
more complete discussion see \cite{Bj3}. Let 
$\mu: L \to \Z$ be the M\"obius function of $L$.

\begin{theorem} For any linear order on $[n]$, 
\begin{enumerate}
\item $\nc(G)$ is a pure simplicial complex of 
dimension $\rk(G)-1$.
\item The cardinality of $\nc^p(G)$ is equal to $w_p(L)$, the $p^{\rm 
th}$ Whitney number of $L$. 
\item For every flat $X$ of $G$, the cardinality of the set 
$\nc_X(G)$ is equal to $(-1)^{\rk(X)}\mu(X)$
\end{enumerate}
\label{purity}
\end{theorem}
 
The relevance of \nbc\ sets to Orlik-Solomon algebras was established by 
several 
authors independently -- see \cite[Section 7.11, \mbox{\S} 7.10]{Bj3}.

\begin{theorem} The set $\{a_S \ | \ S\in \nc(G)\}$ forms a basis for 
$A(G)$.
\label{nbcbase}
\end{theorem}

\medskip
There are several natural ways in which one might attempt to relate 
$\bar{A}(G)$ 
directly to $\bar{G}$, motivated by the various connections between $A(G)$ and 
$G$. (For instance, independent sets in $G$ correspond to nonzero monomials in 
$A(G)$.) None of these seem to work; the difficulties can all be traced back 
to 
the failure of the exchange axiom.   
There is at least an indirect connection between $\bar{A}(G)$ and $\bar{G}$ 
obtained by generalizing  
the following well-known property of \nbc\ sets. Let us impose the natural 
linear order on $[n]$, unless otherwise noted.

\begin{theorem}[{\rm \cite{Bj3}}] An increasing subset 
$S=\{i_1,\ldots,i_p\}\subseteq [n]$ is 
nbc if and only if $i_k=\min \cl(\{i_k,\ldots,i_p\})$ for each 
$1\leq 
k\leq p$.
\label{bjorner}
\end{theorem}

Replacing matroid closure with line-closure, we propose the following.

\begin{definition} An increasing subset $S=\{i_1,\ldots,i_p\}\subseteq [n]$ 
is \nbb\ if and only if $$i_k=\min \lc(\{i_k,\ldots,i_p\})$$
or each $1\leq 
k\leq p$. 
\label{nbb}
\end{definition}

The collection of \nbb\ sets of $G$ will be denoted by $\nb(G)$. Of 
course $\nb(G)$ is dependent only on $\bar{G}$, rather than $G$. 

\begin{theorem} $\nb(G)$ is a simplicial complex, containing 
$\nc(G)$ as a subcomplex.
\label{ncinnb}
\end{theorem}

\begin{proof} The first assertion follows from the monotonicity of 
the 
line-closure operator. The second is a consequence of the fact that 
$\lc(S)\subseteq \cl(S)$ for any subset $S$ of $[n]$.
\end{proof}
We will customarily specify $\nb(G)$ and $\nc(G)$ by listing the 
facets, or maximal simplices.

Because of the lack of exchange, $\nb(G)$ depends heavily on 
the linear ordering of the points.

\begin{example} Let $G$ be the matroid of Example \ref{wheel}. Then, 
with the natural linear order on $[6]$, the facets of $\nb(G)$ are
$$1246,136,135,125,134, \ \text{and} \ 124.$$
If we adopt the linear ordering $2<1<3<4<5<6$, the new \nbb\ complex 
has facets $$246,236,216,235,215,234,\ \text{and} \ 214.$$ In fact,
for this second linear ordering. $\nb(G)=\nc(G)$. 
\qed
\end{example}

We see from this example that the number of \nbb\ sets of a 
fixed size $p$ is not independent of the linear ordering, and the complex 
$\nb(G)$ may fail to be pure. Compare with  Theorem~\ref{purity}(i) and (ii).
We also see that $\nb(G)$
may agree with $\nc(G)$ even when $G$ is not line-closed. However, 
these \nbb\ sets do capture the lack of line-closure, in the following sense.

\begin{theorem} The matroid $G$ is line-closed if and only if 
$\nb(G)=\nc(G)$ for every linear ordering of $[n]$
\label{nb=nc}
\end{theorem}

\begin{proof} Suppose $G$ is not line-closed. Then there exists 
a line-closed set $X$ which is not closed. Let 
$i\in\cl(X)-X$, and choose a linear order on $[n]$ such that $i$ 
precedes $\min(X)$. Now, let $S=(i_1,\ldots,i_p)$ be the 
lexicographically first 
ordered basis for the flat $\cl(X)$ which is contained in $X$. Then, 
by the 
choice of ordering, $S\not \in \nc(G)$, by Theorem \ref{bjorner}. We 
claim $S\in 
\nb(G)$. Suppose not. Then, for some $k$, $i:= \min \lc\{i_k, 
\ldots,i_p\}$ is 
less than $i_k$. Since $X$ is line-closed, $i\in X$. Also, by the 
exchange axiom 
in $G$, $S-\{i_k\}\cup\{i\}$ is a basis for $\cl(X)$, and is 
lexicographically 
smaller than $S$. This contradicts the choice of $S$. Thus 
$S\in\nb(G)-\nc(G)$, 
so $\nc(G)\not = \nb(G)$. Conversely, if $G$ is line-closed, then 
$\nb(G)=\nc(G)$ by Theorem \ref{bjorner}.
\end{proof}
\end{paragraph}

\begin{paragraph}{Independence of \nbb\ monomials in $\bar{A}(G)$}
We turn now to the analysis of the quadratic closure $\bar{A}(G)$ 
of 
the Orlik-Solomon algebra. 

For each line-closed set $X\in \bar{L}(G)$, let $\E_X$ be the 
subspace of 
$\E$ spanned by monomials $e_S$ for which 
$\lc(S)=X$.  Then we have a grading of $\E$ by $\bar{L}(G)$: 
$$\E=\oplus_{X\in \bar{L}(G)} \E_X.$$ Let $\J_X=\J\cap \E_X$ and 
$\bar{A}_X(G)=\E_X/\J_X$. Then we have the following analogue of 
\cite[Theorem 3.26]{OT}.

\begin{lemma} $\bar{A}(G)=\oplus_{X\in\bar{L}(G)} \bar{A}_X(G)$.
\label{sum}
\end{lemma}

\begin{proof} The ideal \J\ is generated by elements $\partial 
e_{ijk}$ where $\{i,j,k\}$ is dependent.  Since $G$ has no multiple 
points, $\{i,j,k\}$ is a circuit. Then 
$\lc(\{i,j\})=\lc(\{i,k\})=\lc(\{j,k\})$, each being equal to 
$\cl(\{i,j,k\})$. This shows that $\partial e_{ijk}$ is 
homogeneous in the grading above. Thus $\J=\oplus_{X\in\bar{L}(G)} 
\J_X$, 
and the result follows.
\end{proof}

We will also use the following elementary observation, proof left to 
the reader.
\begin{lemma} The graded derivation $\partial: \E \to \E$ induces a 
graded derivation $\bar{\partial}: \bar{A}(G) \to \bar{A}(G)$.
\end{lemma}

We are now prepared to prove the main result.
For $S\subset [n]$ we denote by $\bar{a}_S$ 
the image of $e_S$ in the quadratic closure $\bar{A}(G)$.

\begin{theorem} The set $\{\bar{a}_S \ | \ S\in \nb(G)\}$ is linearly 
independent in $\bar{A}(G)$.
\label{indep}
\end{theorem}

\begin{proof}
With Lemma \ref{sum} in hand the proof is identical to the argument 
in 
the proof of Theorem \ref{nbcbase}.  It is enough to prove the result 
for \nbb\ sets of a fixed size $p$.  Then we induct on $p$.  Suppose 
$$\sum_{S\in \nb^p(G)} \lambda_S \bar{a}_S=0.$$
By Lemma \ref{sum} we may assume that 
$\lc(S)=X$ for a fixed element $X\in \bar{L}(G)$ and all $S$ in the 
sum.  
Setting $i_0=\min(X)$, we have $\min(S)=i_0$ for all $S$, by 
definition of $\nb(G)$.  Write $S'=S-\{i_0\}$.  Then we have 
$$\bar{a}_{i_0}\wedge 
\biggl(\underset{\lc(S)=X}{\sum_{S\in\nb^p(G)}} 
\lambda_S 
\bar{a}_{S'}\biggr)=0.$$ Applying the derivation $\bar{\partial}$ we obtain 
$$\underset{\lc(S)=X}{\sum_{S\in\nb^p(G)}} \lambda_S \bar{a}_{S'} + 
\underset{\lc(S)=X}{\sum_{S\in\nb^p(G)}} \bar{a}_{i_0} \wedge 
\bar{\partial}\bar{a}_{S'}=0.$$ Using again the definition of 
$\nb(G)$, we have that $i_0\not\in \lc(S')$ for $\lc(S)=X$.  Then, 
applying Lemma \ref{sum} once more, we have 
$$\underset{\lc(S)=X}{\sum_{S\in\nb^p(G)}}  \lambda_S \bar{a}_{S'} =0.$$ Since 
$S\in \nb^p(G)$ implies 
$S'\in 
\nb^{p-1}(G)$, we conclude $\lambda_S=0$ for 
all $S$ by the inductive hypothesis.
\end{proof}

As a consequence we obtain half of Conjecture \ref{boston}. 
\begin{corollary} Suppose $A(G)$ is quadratic. Then $G$ is 
line-closed.
\label{oneway}
\end{corollary}

\begin{proof} If $G$ is not line-closed, then by Theorems 
\ref{ncinnb} 
and \ref{nb=nc} 
there is a linear ordering of $[n]$ such that the cardinality of 
$\nb(G)$ is strictly greater than that of $\nc(G)$. Then, by 
Theorems \ref{indep} and \ref{nbcbase}, we have $\dim \bar{A}(G)>\dim 
A(G)$, so $A(G)$ is not quadratic.
\end{proof}

Because of Example \ref{wheel}, it is not the case 
that $\{\bar{a}_S \ | \ S\in \nb(G)\}$ spans $\bar{A}(G)$ for every 
linear order. 
When we announced Theorem~\ref{indep} in the  Boston lecture\cite{F13} we 
expressed 
some hope that 
one could show the existence of some linear order for which the \nbb\ 
monomials span $\bar{A}(G)$, yielding a proof of the converse 
of Corollary~\ref{oneway} as well as a combinatorial calculation of 
$\phi_3(G)$ via Corollary~\ref{count}. Subsequently S.~Yuzvinsky 
found a counterexample. 

\begin{example}[{\rm Yuzvinsky \cite{DY}}] Let $G$ be the rank-three matroid 
on 
$[8]$  
with nontrivial lines $$123,\ 148,\ 257,\ 3678,\ \text{and}\ 456,$$ 
pictured in Figure~\vref{braced}.

\begin{figure}

\begin{center}
\epsfig{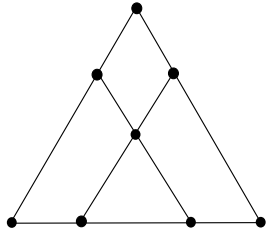}
\caption{$G$ line-closed, $A(G)$ not quadratic}
\label{braced}
\end{center}
\end{figure} 
One can use Theorem~\ref{baseclosure} of Section~\ref{formal} below
to check fairly quickly by hand that 
$G$ 
is line-closed. On the other hand, one computes $\phi_3(G)=16$. (We use a {\em 
Mathematica} script available from the author.) Then, by  
Theorem~\ref{count}, we have $\dim \bar{A}\,^3 = 16$. But $\dim 
A^3(G) = 14$.
Thus $A(G)$ is not quadratic.
\label{yuzex}
\end{example}

For us, Example \ref{yuzex} shows that 
there may be no ``good'' linear order, for which $\nb(G)$ spans 
$\bar{A}(G)$, for some matroids $G$. See \cite{DY} for further 
discussion of 
the converse of Corollary \ref{oneway}. 

\medskip
Using Theorem~\ref{oneway} we have an alternate proof of 
\cite[Prop. 5.1]{F4}.

\begin{corollary} If $n>\rk(G)$ and $G$ has a basis each of whose 
two-point subsets is closed in $G$, then $A(G)$ is not quadratic.
\label{doublepoint}
\end{corollary}

\begin{proof} Such a basis $B$ would form a line-closed set by 
hypothesis, but it cannot be closed in $G$ since 
$|\cl(B)|=n>\rk(G)=|B|$.
\end{proof}

We may also use the work on quadratic algebras together with 
Corollary~\ref{oneway} to 
give a nice sufficient condition for line-closure of matroids. The following 
assertion is a generalization of \cite[Proposition 3.2]{F8}, with essentially 
the same proof.

\begin{corollary} Suppose, for every circuit $S$ of $G$ with 
$|S|\geq 4$,   
the closures in $G$ of two disjoint two-point subsets of $S$ meet.  
Then $G$ is line-closed.
\label{parallel}
\end{corollary}

\begin{proof} We show that $A(G)$ is quadratic by verifying directly that 
$\I^p\subseteq \J^p$ for all $p\geq 3$. 
Let $S=\{i_1,\ldots, i_p\}$ be a circuit, $p\geq 4$, and suppose $i_0\in 
\cl(\{i_1,i_2\})\cap \cl(\{i_3,i_4\})$. Then $$\partial 
e_{i_1i_2i_3i_4}=(e_{i_3}-e_{i_4})\partial e_{i_0i_1i_2} + 
(e_{i_1}-e_{i_2})\partial e_{i_0i_3i_4}.$$ Thus $\partial e_{i_1i_2i_3i_4}$ 
lies in $\J$. Then $e_{i_1i_2i_3i_4} \in \J$, and it follows from the Leibniz 
rule that $\partial e_S$ lies in $\J$.
\end{proof}

If $G$ has rank three, the hypothesis of Theorem \ref{parallel} can be weakened, 
for in this case it suffices to show that $\partial e_S\in \J$ for those 
circuits $S$ with $|S|\geq 4$ and $1\in S$. Arrangements of rank three whose 
matroids satisfy the weaker hypothesis
are called {\em parallel arrangements.} See \cite{FR1}.
\end{paragraph}

\begin{paragraph}{Generalization to high rank/degree}
All of the results of this section on line-closure and quadraticity 
can be generalized, with essentially identical proofs. 

\begin{definition} A subset of $[n]$ is {\em $r$-closed} if it 
contains the closures of all of its $p$-subsets for all $p\leq r$. 
The matroid $G$ is 
{\em $r$-closed} if every $r$-closed set is closed.
\end{definition}

For arbitrary $G$ the collection $L_r(G)$ of $r$-closed sets forms a 
lattice, and we have a sequence of surjective order-preserving maps 
$$B_n=L_1(G)\to \bar{L}(G)=L_2(G) \to \cdots \to L_r(G) \to \cdots 
\to 
L_n(G)=L(G),$$ 
where $B_n$ is the boolean lattice.
 
\begin{definition} The {\em degree $r$ closure} $\bar{A}_r(G)$ of 
$A(G)$ 
is $\E/\J_r$, where $\J_r$ is the ideal generated by the elements of 
\I\ of degree less than or equal to $r$.
\end{definition}

Thus we have a sequence of surjective homomorphisms $$\E=A_1(G) \to 
\bar{A}(G)=A_2(G) \to \cdots \to A_r(G) \to \cdots A_n(G)=A(G).$$

\begin{definition} An ordered subset $S=\{i_1,\ldots, i_p\}$ is 
$r$-\nbb\ if
for each $k$, $i_k$ is the first element in the $r$-closure of 
$\{i_k\ldots, i_p\}$.
\end{definition}

\begin{lemma} $A_r(G)=\oplus_{X\in\bar{L_r}(G)} \bar{A}_X(G)$.
\end{lemma}

\begin{proof} The crucial points are (i) that $\J_r$ is generated by 
boundaries of {\em circuits} of size at most $r+1$, and (ii) that 
$r$-closure agrees with matroid closure on sets of size at most $r$.
Using these observations the proof of Lemma~\ref{sum} is easy to adapt 
to the more general setting.
\end{proof}

The proof of the following generalization is now identical to the proof of 
Theorem~\ref{indep}.
\begin{theorem} The set of monomials in $\A_r(G)$ corresponding to 
the $r$-\nbb\ sets of $G$ forms a linearly independent set.
\end{theorem}

\begin{corollary} If $A_r(G)=A(G)$, then $G$ is $r$-closed.
\end{corollary}
\end{paragraph}

\begin{paragraph}{The \nbb\ complex and a generalization to nonlinear orders}
After formulating Definition \ref{nbb} and proving 
Theorem~\ref{indep}, we found that our notion of \nbb\ set
coincides with a special case of the more
general notion of \NBB\ set in a finite lattice with a partial 
ordering of the atoms, introduced by Blass and Sagan in \cite{BSag}.
These results are stated only for the line-closure of $G$, but 
again analogous results will hold for $r$-closure.

\begin{definition}(\cite{BSag}) Suppose $(\bar{L},\leq)$ is a finite lattice, 
and $\preceq$ is a partial ordering of the atoms of $\bar{L}$. A set $T$ of 
atoms is {\em bounded below} if there exists an atom $a$ such that $a< 
\bigvee T$ 
and $a\prec t$ for all $t\in T$. A set $S$ is \NBB\ if $S$ does not contain 
any set $T$ which is bounded below.
\end{definition}

\begin{theorem} A set $S$ is \nbb\ if and only if $S$ is an \NBB\ set 
in the lattice $\bar{L}(G)$ for the given linear ordering of the atoms.
\label{sagan}
\end{theorem}

\begin{proof} In our setting the atom ordering $\preceq$ is the natural linear 
ordering on $[n]$. In this context a set $T\subseteq [n]$ is bounded below if 
and only if there exists $i\in 
\lc(T)$ with $i<\min(T)$.  Suppose $S=\{i_1,\ldots,i_p\}$ is \nbb\, 
and 
$T\subseteq S$.  Let $i_k=\min(T)$.  Then $\lc(T)\subseteq 
\lc(\{i_k,\ldots,i_p\}$, so $\min(\lc(T))\geq 
\min(\lc(\{i_k,\ldots,i_p\}))$.  Since $S$ is \nbb\, we conclude 
$\min(\lc(T))=i_k$, so $T$ is not bounded below.  Conversely, if $S$ 
is not \nbb\, then for some $k$, $T=\{i_k,\ldots,i_p\}$ is bounded 
below by $\min(\lc(T))<i_k,$ and thus $S$ is not \NBB.
\end{proof}

This observation yields a numerical result on the number of \nbb\ 
sets, by one of the main results of \cite{BSag}.
 Let $\bar{\mu}: \bar{L}(G) \to \Z$ denote 
the 
M\"obius function of $\bar{L}(G)$.

\begin{corollary} The sum $\underset{\lc(S)=X}{\underset{S\in \nb(G)}{\sum}} 
(-1)^{|S|}$ is equal to 
$\bar{\mu}(X)$. 
\end{corollary}

Let $\NB(G,\preceq)$ denote the collection of \NBB\ 
sets of $\bar{L}(G)$ under the atom-order $\preceq$.

\begin{theorem} Suppose $\preceq$ is a partial order on $[n]$ with the 
property that each line-closed set $X$ has a unique smallest element 
relative to $\preceq$. Then $$\{\bar{a}_S \ | \ S\in 
\NB(G,\preceq)\}$$ is linearly independent in $\bar{A}(G)$.
\label{nonlinear}
\end{theorem}

\begin{proof} Assume without loss that the natural order on $[n]$ is a 
linear extension of $\preceq$. Then the proof of \ref{indep} goes 
through without change.
\end{proof}

We also have the following analogue of Theorem~\ref{nb=nc}. Recall 
that $\nc(G)$ is determined by a linear order on $[n]$.

\begin{theorem} $\nc(G)\subseteq\NB(G,\preceq)$ for any linear order 
which extends $\preceq$.
\end{theorem}

\begin{proof}
Suppose $T$ is a bounded below set. Let $a\leq \bigvee T$ with 
$a\prec t$ for all $t\in T$. Then $a$ precedes $\min(T)$ in any linear 
extension of $\preceq$. Since $\bigvee T =\lc(T)\subseteq \cl(T)$, it 
follows that $T$ contains a broken circuit.
\end{proof}

Theorem~\ref{nonlinear} raises the possibility of finding more than 
$|\nc(G)|$ linearly independent monomials in $\bar{A}(G)$, even when  
$G$ is line-closed. 
At this point we have no examples of this phenomenon.
 
\end{paragraph}
\end{section}

\begin{section}{Combinatorial notions of formality}
\label{comb}

Our research on line-closed matroids was motivated by the study of 
formal arrangements, and specifically by attempts to describe 
combinatorially the property of an arrangement being ``generic with 
given codimension two structure." We start this section by recalling the 
definition of formal arrangement, and outlining some of the motivation and 
main results. We then turn to various combinatorial versions of 
formality. These have been studied to some extent before, but the definitions 
have never been published. We present some newly rediscovered results 
and examples, which appeared long ago in the combinatorial literature but not 
in 
the context of formal arrangements. 

Line-closure turns out to be the 
strongest among the properties we study here. 
The remaining notions form a hierarchy 
descending to formality of a realization, which in itself is not 
a combinatorial notion. We show by example that none of the combinatorial 
notions have the nice topological or algebraic consequences 
that formality affords. 

\begin{paragraph}{Formal arrangements}
The notion of formality was introduced in \cite{FR1} and studied further in 
\cite{Y1,BT,BB}. The terminology is unfortunate; there 
is a notion of formality of spaces that is important in rational homotopy 
theory, and has implications for arrangements, but the definition of formal 
arrangement is completely unrelated.

We adopt the definition from \cite{BT}. Henceforth 
assume \A\ is an {\em essential} arrangement, that is, 
$\rk(G)=\dim(V)$. Let $e_i$ 
denote the $i^{\rm th}$ standard basis vector in $\K^n$. The {\em 
weight} of a vector in $\K^n$ is the number of nonzero entries.

\begin{definition}
Let $\rho: \K^n \to V^*$ be the linear mapping defined by 
$\rho(e_i)=\al_i$. Then \A\ is {\em formal} if the kernel of $\rho$ 
is 
spanned by elements of weight at most three.
\label{formaldef}
\end{definition}

An element of $\ker(\rho)$ of weight three corresponds to a dependent 
set of $G$ of size three, and thus of rank two.  A vector in $\ker(\rho)$ 
gives the coefficients in a 
dependence relation among the linear forms. Thus \A\ is 
formal if all dependence relations among the forms $\al_i$ are 
consequences of ``rank two dependence relations." 

There is a natural way to associate a subspace $W$ of $\K^n$ of 
dimension 
$r$ with a (possibly degenerate) arrangement of $n$ hyperplanes in 
$\K^r$, by 
considering 
the set of intersections of $W$ with the $n$ coordinate hyperplanes 
as an 
arrangement in $W$. 
If $K\subseteq \K^n$ denotes the kernel of $\rho$, then its 
orthogonal complement $K^\perp$ returns the original arrangement \A\ 
under this construction. Indeed, the linear mapping 
$$\Phi=(\al_1,\ldots,\al_n): V \to \K^n$$ carries $V$ isomorphically 
to 
$K^\perp$ and $H_i$ to the intersection of $K^\perp$ with 
$\{x_i=0\}$. Let $F\subseteq K$ denote the subspace spanned by 
elements of weight three. The arrangement corresponding to the 
subspace $F^\perp \subseteq \K^n$ is called the {\em formalization} 
of \A, denoted $\A_F$. This construction first appears in 
\cite{Y1}.

A {\em section} $\bar{\B}$ of an arrangement \B\ in $V$ is 
formed by 
intersecting the 
hyperplanes of \B\ with a linear subspace $W$ of $V$. The 
section $\bar{\B}$ is {\em generic} if $W$ is transverse to every 
intersection of hyperplanes of \B. In this case
the combinatorics and topology of $\bar{\B}$ depend only on \B\ and
$\dim(W)$. A section of \B\ by a 3-dimensional subspace is  
called a 
{\em planar section}. 
 
\begin{theorem}[{\rm \cite{Y1}}] Let \A\ be an essential arrangement. Then
\begin{quote}
\begin{enumerate}
\item $\A_F$ is formal.
\item $\A$ is a section of $\A_F$.
\item $\A$ and $\A_F$ have identical generic planar sections.
\item $\A$ is formal if and only if $\rk(\A)=\rk(\A_F)$.
\end{enumerate}
\end{quote}
\label{formalprop}
\end{theorem}

We remark that \A\ need not be a generic section of $\A_F$, nor of any other 
arrangement. Indeed, if there are no nontrivial 
lines in $G$, then $\A_F$ is the boolean arrangement. Then, if \A\ has some 
nontrivial plane, and has more than four elements, \A\ will not be a generic 
section of $\A_F.$ If \A\ is inerectible, it will not be a generic section of 
any arrangement. Such arrangements are easy to construct.

The interest in formal arrangements is due to the following theorem.

\begin{theorem} Let \A\ be an essential arrangement. Then
\begin{quote}
\begin{enumerate}
\item  If \A\ is a $K(\pi,1)$ arrangement, then \A\ is 
formal.
\item  If \A\ is a rational $K(\pi,1)$ arrangement, then \A\ 
is 
formal.
\item  If \A\ is a free arrangement, then \A\ is formal.
\item  If \A\ has quadratic Orlik-Solomon algebra, then \A\ is formal.
\end{enumerate}
\end{quote}
\label{goodstuff}
\end{theorem}
Assertions (i) and (iv) above are easy consequences of 
Theorem~\ref{formalprop}, 
and (ii) is a consequence of (iv), because rational $K(\pi,1)$ arrangements 
have 
quadratic Orlik-Solomon algebras. See \cite{FR2}. Assertion (iii) was proved 
in \cite{Y1}.

In \cite{Y1}, Yuzvinsky presented examples of two arrangements, one formal and 
the other not, with the same underlying matroid. See \cite{FR4} for 
diagrams of the dual point configurations. 
The underlying matroid is the 
dual of the matroid of complete bipartite graph $K_{3,3}$. This observation 
yields different realizations and a geometric explanation of this phenomenon.

\begin{example}
\label{k33}
The diagram on the left in Figure~\vref{yuzfig} is a representation of the 
rank-four matroid $G^*$ dual to the graphic matroid of $K_{3,3}$. By 
considering 
intersecting planes (in $\C^3$), one can see that the three dotted lines in 
this, or in any \C -representation of $G^*$, must be concurrent. The diagram 
on 
the 
right is a representation of the truncation $T(G^*)$ in which the 
corresponding 
lines are not concurrent. And indeed, one can show that the configuration on 
the 
right is formal. In fact, it cannot be lifted to a rank four configuration 
with the same points and lines.


\begin{figure}[h]
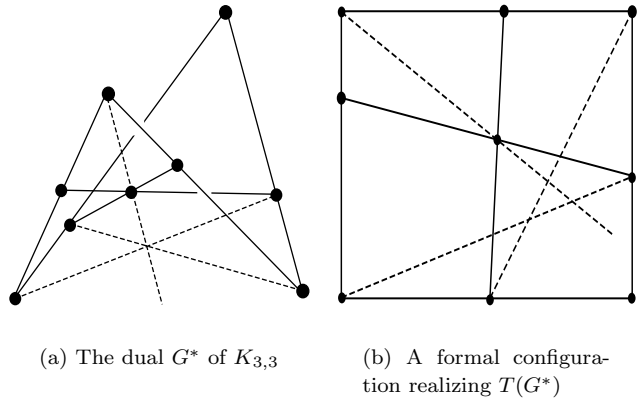

\begin{center}
\mbox{\subfigure[The dual $G^*$ of 
$K_{3,3}$]{\epsfig{file=k33.epsf,height=4cm, width=4cm}}\quad
\subfigure[A formal configuration realizing 
$T(G^*)$]{\epsfig{file=nonformal.epsf,height=4cm,width=4cm}}}
\caption{Formality is not matroidal}
\label{yuzfig}
\end{center}
\end{figure}

\end{example}

\end{paragraph}

\begin{paragraph}{Combinatorial formality}
\label{formal}
Example \ref{k33} shows that formality is not a 
combinatorial property. Since the discovery of these examples, 
efforts have been made to strengthen Theorem~\ref{goodstuff} by replacing the 
formality  assumption with some purely combinatorial property. In this 
subsection 
we present several reasonable candidates. 

\medskip
Theorem~\ref{formalprop}(ii) suggests a natural combinatorial formulation of 
formality.  We recall for the reader the notion of strong map (or 
quotient) of matroids. See \cite[Section 7.4 and Chaps. 8-9]{Wh1} for 
more details. 
Suppose $G'$ and $G$ are two matroids on ground set $[n]$.  We say 
$G$ is a {\em quotient} of $G'$ if every closed set in $G$ is 
closed in $G'$. This is the case precisely when the identity map on 
$[n]$ is a strong map from $G'$ to $G$. If $G'$ is the matroid of an 
arrangement \A, 
then the matroid of any section of \A\ is a quotient of $G$.  The 
matroid of a generic $r$-dimensional section of \A\ coincides with the 
{\em truncation} $G^{[r]}$ of the matroid $G$ to rank $r$, the matroid 
whose dependent subsets are those of $G$ together with every subset 
of 
size greater than $r$.

\begin{definition} A matroid $G$ is {\em taut} if $G$ is not a 
quotient of any matroid $G'\not=G$ satisfying $(G')^{[3]}=G^{[3]}$. 
\end{definition}

The condition $(G')^{[3]}=G^{[3]}$ says merely that $G$ and $G'$ have the 
same points and lines.
The following assertion is a consequence of Theorem 
\ref{formalprop}.

\begin{corollary} If $G$ is a taut matroid, then every arrangement 
realizing $G$ is formal.
\label{cformal}
\end{corollary}

In a lecture in 1992 \cite{Y6}, Yuzvinsky formulated a definition of 
``dimension" 
of a geometric lattice, or matroid, based on line-closure: the {\em dimension} 
of $G$ is the size of the smallest set of points whose line-closure is $[n]$.
Corollary~\ref{look} below was presented in \cite{Y6}, but was never 
published.
The result was already known to matroid theorists \cite{Cra}.

\begin{theorem} Suppose $G$  has a basis (of $\rk(G)$ points) whose 
line-closure 
is $[n]$. Then $G$ is taut.
\label{tautness}
\end{theorem}

\begin{proof} Suppose $G$ is a quotient of a matroid $G'$ with the 
same points and lines as $G$.  Since the closure of $B$ in $G'$ 
contains the 
line-closure of $B$ in $G'$, which agrees with the line-closure of 
$B$ 
in $G$, we conclude that $B$ is a basis for $G'$.  Thus 
$\rk(G')=\rk(G)$. It follows that $G'=G$.
\end{proof}

The following criterion is the standard method to prove an arrangement is 
formal, 
although it has never appeared in the literature.

\begin{corollary} Suppose $G$ has a basis whose line-closure is $[n]$. Then 
every realization of $G$ is formal.
\label{look}
\end{corollary}

\medskip
We now have four notions which might capture the combinatorics of 
formality, at least in spirit:

\begin{enumerate}
\item $G$ is line-closed,
\item $G$ has a basis whose line-closure is $[n]$,
\item $G$ is taut, that is, $G$ is not a proper quotient preserving points and 
lines, and
\item every realization of $G$ is formal.
\end{enumerate}

We have the string of implications $$(i)\implies (ii) \implies (iii) \implies 
(iv).$$

The first of these is trivial, and the others were proved in the preceding 
paragraphs. We now give counter-examples for the first two of the reverse 
implications.
We note that Example~\ref{k33} provides a formal arrangement whose matroid 
fails 
to satisfy (iv).

\begin{example} In \cite{Cra} there appears an example of a 
matroid $G$ of rank three on nine points, with the property that 
no set of three points line-closes to the entire ground set $[9]$. 
See Figure~\vref{crapo}. Thus (ii) fails. But it is easy to 
see that $G$ is not a quotient of a rank-four matroid with the same 
points and lines, so (iii) is satisfied.


\begin{figure}[h]
\begin{center}
\epsfig{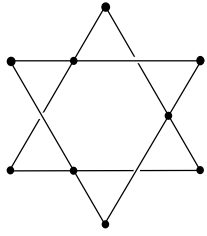}
\caption{Crapo's matroid}
\label{crapo}
\end{center}
\end{figure}

\label{threenottwo}
\end{example}

\begin{example} The rank-three wheel, illustrated in 
Figure~\vref{wheel}, provides an example of a matroid which satifies 
(ii) but is not line-closed.
\label{twonotone}
\end{example}

Finding a counter-example for the implication ``$(iv) \implies (iii)$'' 
presents a delicate problem. One needs a \C-representable matroid 
which is a proper quotient, preserving points and lines, but such that 
the quotient ``mapping'' is not representable over \C, either because 
the larger matroid is not \C-representable, or because the  larger 
matroid is not realizable in such a way that the original matroid is 
obtained from it via projection.

\end{paragraph}

\begin{paragraph}{Local formality}
In \cite{Y1} an arrangement is defined to be {\em locally formal} if,
for every flat $X$, the arrangement $\A_X:=\{H_i \ | \ i \in
X\}$ is formal. This idea can be applied to the 
combinatorial notions of formality discussed above. We will focus on 
the local version of tautness, for reasons that will become clear.

\begin{definition} A matroid $G$ is {\em
locally taut} if the restriction of $G$ to any flat is taut.
\end{definition}

By Corollary~\ref{cformal}, any 
realization of a locally taut matroid is locally formal. 
Theorem~\ref{tautness} can be used to compare local tautness to 
line-closure, in the next pair of results.

\begin{corollary} Suppose $G$ is a matroid in which every flat $X$ 
has a basis whose line-closure is equal to $X$. Then $G$ is locally 
taut.
\label{localtaut}
\end{corollary}

\begin{proof} Thus is an immediate consequence of 
Theorem~\ref{tautness}.
\end{proof}

By way of contrast, the definition of line-closed matroid may be 
restated as follows.

\begin{theorem}
A matroid $G$ is line-closed if and only if, for every flat $X$, 
{\em every} 
basis of $X$ has line-closure equal to $X$.
\label{baseclosure}
\end{theorem}

Thus every line-closed matroid is even locally taut, strengthening 
Theorem~\ref{tautness}. The rank-three wheel $W_3$ (Examples~\ref{wheel} and 
\ref{twonotone}) serves as an example of a locally taut matroid which 
is not line-closed.

\medskip
We make one more observation concerning locally taut matroids.

\begin{theorem} Suppose $G$ is locally taut. Then $G$ is determined 
by its points and lines.
\label{pointline}
\end{theorem}

\begin{proof} This is a consequence of part (5) of Theorem 7.5.4 in 
\cite{Bry3}, which asserts that $G$ is determined by its essential 
flats, along with their ranks. A flat of $G$ is essential if it 
is a truncation of a matroid of higher rank. Truncation is in 
particular a quotient map, and preserves points and lines, so long as 
the image has rank at least two. Thus, if $G$ is locally taut, 
the only essential flats of $G$ are the nontrivial lines. So $G$ is 
determined by the number of points and the list of nontrivial lines.
\end{proof}

The converse of this result also holds. For if $G$ is not locally 
taut, then $G$ is not taut, so $G$ is a quotient of a distinct 
matroid with the same points and lines. A natural question is 
evident: what axioms govern the point-line incidence structures of 
locally taut matroids?

\end{paragraph}
\begin{paragraph}{Topology vs. combinatorics}

One motivation of the present discussion is to find combinatorial 
conditions with the same consequences as formality, as in 
Theorem~\ref{goodstuff}. Unfortunately, we have no positive results 
of this sort, aside from the implication ``quadratic Orlik-Solomon 
algebra implies line-closed matroid'' strengthening assertions (ii) and 
(iv) of that theorem.

First of all, it is not the case that $K(\pi,1)$ 
or free arrangements must have line-closed matroids.

\begin{example} Consider the arrangement \A\ with defining equation 
$Q(x,y,z)=z(x+y)(x-y)(x+z)(x-z)(y+z)(y-z)$, denoted $J_2$ in 
\cite[section 2.6]{FR1}. The underlying matroid is the non-Fano 
plane, pictured in Figure~\vref{nonfano}. Then \A\ is 
not line-closed: apply Corollary~\ref{doublepoint} to the three 
``edge-midpoints.'' But this arrangement is well-known to be both 
free and $K(\pi,1)$ - see \cite{FR1}.


\begin{figure}[h]
\begin{center}
\epsfig{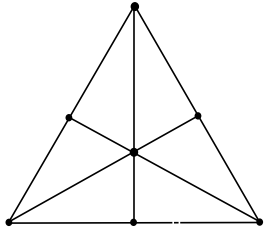}
\caption{Not line-closed, with a free, $K(\pi,1)$ realization}
\label{nonfano}
\end{center}
\end{figure}

\label{jay2}
\end{example}

Line-closure is not sufficient for $K(\pi,1)$-ness or freeness either.
\begin{example} Let \A\ be the parallel arrangement with defining equation 
$Q(x,y,z)=
z(x+z)(x-z)(y+z)(y-z)(x+y+2z)(x+y-2z)$, which appears as $X_2$ in 
\cite[Section 2.6]{FR1}. 
The underlying matroid $G$ of \A\ is pictured in Figure~\vref{kohno}. 


\begin{figure}[h]
\begin{center}
\epsfig{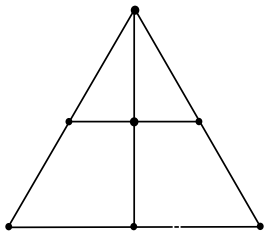}
\caption{Line-closed, with a non-free, non-$K(\pi,1)$ realization}
\label{kohno}
\end{center}
\end{figure}


Since the hypothesis of Theorem \ref{parallel} is 
satisfied, 
$G$ is line-closed. But  
\A\ is not a free arrangement, nor is \A\ a $K(\pi,1)$ arrangement, nor a 
rational $K(\pi,1)$ arrangement (though it has quadratic 
Orlik-Solomon algebra). The first of these assertions holds because the 
characteristic polynomial of $G$ has non-integer roots. The second assertion 
was 
proved by L.~Paris (unpublished), and the last was proved in \cite{Y5}. See 
\cite{FR4}.
\end{example}

We close with a fascinating conjecture about the combinatorial structure of 
free or $K(\pi,1)$ arrangements inspired by the preceding discussion. 
First we note the stronger version of Theorem~\ref{goodstuff}: 
a free or $K(\pi,1)$ arrangement 
must in fact be {\em locally} formal \cite{F1,Y1}. 

Now, let us consider the hierarchy 
``$(i)\implies (ii)\implies (iii)\implies (iv)$'' among the 
combinatorial notions we have introduced. We have seen that 
that the matroids of $K(\pi,1)$ or free arrangements need not satisfy 
(i). The question whether such matroids satisfy (iv) is very close to 
two famous and important problems in the theory of arrangements. Indeed, 
a matroid which underlies both a free arrangement 
and also a non-formal arrangement would be a counter-example to Terao's 
conjecture, that freeness is a matroidal property. A matroid which 
underlies both a $K(\pi,1)$ arrangement and a non-formal arrangement 
would improve upon Rybnikov's (non-$K(\pi,1)$) counter-examples to the 
homotopy-type conjecture, that the homotopy type of the complement is 
matroidal. See \cite{FR4} for a discussion of the latter problem.

Finally, we note that the matroid appearing in Example~\ref{jay2}, while 
not line-closed, is taut (i.e., satisfies (iii)), in fact locally 
taut, by Theorem~\ref{look}, and is therefore determined by its points and 
lines. 

\begin{conjecture} The underlying matroid of a 
free or 
$K(\pi,1)$ arrangement is determined by its points and lines.
\end{conjecture}

\end{paragraph}
\end{section}

\begin{ack} These ideas have developed through conversations with 
several people over the years. We are especially grateful to Sergey 
Yuzvinsky and Joseph Kung for helpful discussions along the way. In particular, 
Kung pointed out Theorem \ref{pointline} to us. Samantha Melcher helped with 
examples related to the \nbb\ complex during a summer REU project in 1998.  We 
also thank Henry Crapo for explaining to us 
his interpretation of Yuzvinsky's examples, which we reproduced in 
Example \ref{k33}.

\end{ack}

\begin{thebibliography}{Fal89b}

\bibitem[BB97]{BB}
M.~Bayer and K.~Brandt.
\newblock Discriminantal arrangements, fiber polytopes, and formality.
\newblock {\em Journal of Algebraic Combinatorics}, 6:229--246, 1997.

\bibitem[Bj{\"o}92]{Bj3}
A.~Bj{\"o}rner.
\newblock Homology and shellability of matroids and geometric lattices.
\newblock In N.~White, editor, {\em Matroid Applications}, volume~40 of {\em
  Encyclopaedia of Mathematics and Its Applications}, pages 226--283. Cambridge
  University Press, Cambridge, 1992.

\bibitem[Bry86]{Bry3}
T.~Brylawski.
\newblock Constructions.
\newblock In N.~White, editor, {\em Theory of Matroids}, volume~26 of {\em
  Encyclopaedia of Mathematics and Its Applications}, pages 127--223. Cambridge
  University Press, Cambridge, 1986.

\bibitem[BS97]{BSag}
A.~Blass and B.~Sagan.
\newblock M{\"o}bius functions of lattices.
\newblock {\em Advances in Mathematics}, 127:94--123, 1997.

\bibitem[BT94]{BT}
K.~Brandt and H.~Terao.
\newblock Free arrangements and relation spaces.
\newblock {\em Discrete and Computational Geometry}, 12:49--63, 1994.

\bibitem[Cra70]{Cra}
H.~Crapo.
\newblock Erecting geometries.
\newblock In {\em Proceedings of the Second Annual Chapel Hill Conference on
  Combinatorial Mathematics and Its Applications}, pages 74--99, 1970.

\bibitem[DY00]{DY}
G.~Denham and S.~Yuzvinsky.
\newblock Annihilators of ideals in exterior algebras.
\newblock preprint, 2000.

\bibitem[Fal88]{F4}
M.~Falk.
\newblock The minimal model of the complement of an arrangement of hyperplanes.
\newblock {\em Transactions of the American Mathematical Society},
  309:543--556, 1988.

\bibitem[Fal89a]{F8}
M.~Falk.
\newblock The cohomology and fundamental group of a hyperplane complement.
\newblock In {\em Singularities}, volume~90 of {\em Contemporary Mathematics},
  pages 55--72. American Mathematical Society, 1989.

\bibitem[Fal89b]{F6}
M.~Falk.
\newblock On the algebra associated with a geometric lattice.
\newblock {\em Advances in Mathematics}, 80:152--163, 1989.

\bibitem[Fal95]{F1}
M.~Falk.
\newblock {$K(\pi,1)$} arrangements.
\newblock {\em Topology}, 34:141--154, 1995.

\bibitem[Fal99]{F13}
M.~Falk.
\newblock Quadratic algebras and line-closed matroids, June, 1999.
\newblock invited 60-minute lecture, Arrangements in Boston, Northeastern
  University, Boston.

\bibitem[FR85]{FR2}
M.~Falk and R.~Randell.
\newblock The lower central series of a fiber-type arrangement.
\newblock {\em Inventiones mathematicae}, 82:77--88, 1985.

\bibitem[FR86]{FR1}
M.~Falk and R.~Randell.
\newblock On the homotopy theory of arrangements.
\newblock In {\em Complex Analytic Singularities}, volume~8 of {\em Advanced
  Studies in Mathematics}, pages 101--124. North Holland, 1986.

\bibitem[FR00]{FR4}
M.~Falk and R.~Randell.
\newblock On the homotopy theory of arrangements, {II}.
\newblock In M.~Falk and H.~Terao, editors, {\em Arrangements -- Tokyo, 1998},
  Advanced Studies in Mathematics, Tokyo, 2000. Kinokuniya.

\bibitem[OS80]{OS1}
P.~Orlik and L.~Solomon.
\newblock Topology and combinatorics of complements of hyperplanes.
\newblock {\em Inventiones mathematicae}, 56:167--189, 1980.

\bibitem[OT92]{OT}
P.~Orlik and H.~Terao.
\newblock {\em Arrangements of Hyperplanes}.
\newblock Springer Verlag, Berlin Heidelberg New York, 1992.

\bibitem[Oxl92]{Ox}
J.~Oxley.
\newblock {\em Matroid Theory}.
\newblock Oxford University Press, Oxford New York Tokyo, 1992.

\bibitem[PY99]{Y5}
S.~Papadima and S.~Yuzvinsky.
\newblock On rational {$K(\pi,1)$} spaces and {K}oszul algebras.
\newblock {\em Journal of Pure and Applied Algebra}, 144:157--167, 1999.

\bibitem[Rot64]{Rota}
G.-C. Rota.
\newblock On the foundations of combinatorial theory i: Theory of {M}{\"o}bius
  functions.
\newblock {\em Zeitschrift Wahrscheinlichkeitstheorie}, 2:340--368, 1964.

\bibitem[Whi86]{Wh1}
N.~White, editor.
\newblock {\em Theory of Matroids}.
\newblock Cambridge University Press, Cambridge, 1986.

\bibitem[Yuz]{Y6}
S.~Yuzvinsky.
\newblock All arrangements, free and formal.
\newblock Combinatorics Seminar, Mittag-Leffler Institute, Dj{\"u}rsholm,
  Sweden, 1992.

\bibitem[Yuz93]{Y1}
S.~Yuzvinsky.
\newblock First two obstructions to freeness of arrangements.
\newblock {\em Transactions of the American Mathematical Society},
  335:231--244, 1993.

\end{thebibliography}

\bigskip
\makeatletter
{\small \noindent
Department of Mathematics and Statistics\\
Northern Arizona University\\
Flagstaff, AZ  86011-5717 \\
michael.falk@nau.edu
}
\makeatother 

\end{document}